\documentclass[12pt]{article}
\usepackage[psamsfonts]{amssymb}
\usepackage{amsthm}

\begin{document}

\title{\bf On perfect, amicable, and sociable chains}

\author{Jean-Luc Marichal \\
\small{Applied Mathematics Unit, University of Luxembourg} \\
\small{162A, avenue de la Fa\"{\i}encerie, L-1511 Luxembourg, Luxembourg} \\
\small{jean-luc.marichal[at]uni.lu}}

\date{\small{Revised version, September 22, 2000}\footnote{As this paper was
accepted for publication, the author found out that the main problem (Theorem~\ref{thm:list}) was already addressed and solved by Sallows and
Eijkhout \cite{SaEi86}; see also \cite{Ka75,McWa82,SaSh97}.}}

\maketitle

\theoremstyle{plain}
\newtheorem{theorem}{Theorem}
\newtheorem{lemma}[theorem]{Lemma}
\newtheorem{proposition}[theorem]{Proposition}
\newtheorem{corollary}[theorem]{Corollary}

\theoremstyle{definition}
\newtheorem{definition}[theorem]{Definition}
\newtheorem{example}[theorem]{Example}

\theoremstyle{remark}
\newtheorem*{conjecture}{\indent Conjecture}
\newtheorem*{remark}{\indent Remark}

\newcommand{\nit}{\mathbb{N}}                               
\newcommand{\N}{{\cal N}}
\newcommand{\bfx}{\mathbf{x}}
\newcommand{\bfz}{\mathbf{z}}

\begin{abstract}
Let $\bfx = (x_0,\ldots,x_{n-1})$ be an $n$-chain, i.e., an $n$-tuple of non-negative integers $< n$. Consider the operator $s: \bfx \mapsto
\bfx' = (x'_0,\ldots,x'_{n-1})$, where $x'_j$ represents the number of $j$'s appearing among the components of $\bfx$. An $n$-chain $\bfx$ is
said to be perfect if $s(\bfx) = \bfx$. For example, (2,1,2,0,0) is a perfect 5-chain. Analogously to the theory of perfect, amicable, and
sociable numbers, one can define from the operator $s$ the concepts of amicable pair and sociable group of chains. In this paper we give an
exhaustive list of all the perfect, amicable, and sociable chains.
\end{abstract}


\section{Introduction}

Let $n \geqslant 1$ be an integer and let $N := \{0,1,\ldots,n-1\}$. An $n$-{\em chain} is an $n$-tuple
$$
\bfx = (x_0,x_1,\ldots,x_{n-1}),
$$
with $x_i \in N$ for all $i \in N$. Since such an $n$-tuple can be viewed as a mapping from $N$ into itself, the set of all $n$-chains will be
denoted $N^N$, and its cardinality is $|N^N| = n^n$.

Let $2^N$ represent the set of all subsets of $N$. For any $j \in N$, define $S_j : N^N \to 2^N$ as
$$
S_j(\bfx) := \{i \in N \mid x_i = j\}.
$$
Clearly, for any $\bfx \in N^N$, $\{S_j(\bfx) \mid j \in N\}$ is a partition of $N$.

We then say that $\bfx \in N^N$ is a {\em perfect chain} if
$$
x_j = |S_j(\bfx)|, \qquad j \in N.
$$
In other terms, $\bfx \in N^N$ is a perfect chain if, for any $j \in N$, $x_j$ represents the number of $j$'s occuring in
$\{x_0,x_1,\ldots,x_{n-1}\}$. For instance
$$
\bfx = (2,1,2,0,0)
$$
is a perfect 5-chain.

We say that $\bfx, \bfx' \in N^N$ ($\bfx \neq \bfx'$) form a pair of {\em amicable chains} if
$$
\begin{array}{ll}
x'_j = |S_j(\bfx)|, & \quad j \in N, \\
x_j = |S_j(\bfx')|, & \quad j \in N.
\end{array}
$$
For instance
$$
\bfx = (2,3,0,1,0,0) \quad \mbox{and} \quad \bfx' = (3,1,1,1,0,0)
$$
form a pair of amicable 6-chains.

Now, consider the {\em counting operator} $s : N^N \to \{0,1,\ldots,n\}^N$ defined by $\bfx' = s(\bfx)$ with
$$
x'_j = |S_j(\bfx)|, \qquad j \in N.
$$
Given an integer $l \geqslant 3$, we say that the chains $\bfx^{(0)},\bfx^{(1)},\ldots,\bfx^{(l-1)} \in N^N$, satisfying
$$
\bfx^{(k+1)} = s(\bfx^{(k)}),\qquad k \in \{0,\ldots,l-2\},
$$
form a group of $l$ {\em sociable chains} if they are distinct and $s(\bfx^{(l-1)}) = \bfx^{(0)}$. For instance
$$
\bfx^{(0)} = (3,3,0,0,1,0,0), \quad \bfx^{(1)} = (4,1,0,2,0,0,0), \quad \bfx^{(2)} = (4,1,1,0,1,0,0)
$$
form a group of three sociable 7-chains.

Notice that these concepts present some analogies with perfect, amicable, and sociable numbers, see e.g.\ \cite{Cohen70,teRiele83}. Consider the
function $s(n) = \sigma(n)-n$, where $\sigma$ denotes the divisor sum function. A positive integer $n$ is said to be perfect if $s(n) = n$. For
example, 6 is perfect. Two positive integers $m$ and $n$ are said to be amicable if $s(m) = n$ and $s(n) = m$. For example, 220 and 284 are
amicable. An $l$-tuple ($l \geqslant 3$) of positive integers $(n_0,\ldots,n_{l-1})$, satisfying $n_{k+1}=s(n_{k})$ for all $k$, is a sociable
group if these integers are distinct and $s(n_{l-1}) = n_0$. For example, (12~496, 14~288, 15~472, 14~536, 14~264) is a group of 5 sociable
numbers.

The main aim of this paper is to determine all the perfect, amicable, and sociable chains. These are gathered in Theorem~\ref{thm:list} below.
We also investigate the counting operator and point out some of its properties.

The outline of this paper is as follows. In Section 2 we determine conditions under which the iterates of the counting operator are well
defined. In Section 3 the results are presented of an exhaustive computation of all the perfect, amicable, and sociable chains. Finally, Section
4 is devoted to a description of the range of the counting operator and its iterates.

\section{Preliminary results}

In this section we investigate the counting operator $s$ introduced above as well as its iterates. We first observe that this operator does not
always range in $N^N$. For example, if $n = 4$, we have
$$
s(2,2,2,2) = (0,0,4,0) \notin N^N.
$$

We thus need to restrict the domain of $s$ to chains $\bfx$ such that each element of the infinite sequence
$$
\bfx, s(\bfx), s(s(\bfx)), s(s(s(\bfx))),\ldots
$$
belongs to $N^N$. The following results deal with this issue.

\begin{lemma}
Let $\bfx \in N^N$ and $\bfx' = s(\bfx)$. Then
\begin{eqnarray}
\sum_{j \in N} x'_j & = & n, \label{eq:sum/xj} \\
\sum_{j \in N} j \, x'_j & = & \sum_{j \in N} x_j. \label{eq:sum/jxj}
\end{eqnarray}
\end{lemma}

\begin{proof}
Since $\{S_j(\bfx) \mid j \in N\}$ is a partition of $N$, we simply have
$$
\sum_{j \in N} x'_j = \sum_{j \in N} |S_j(\bfx)| = |N| = n ,
$$
and, by counting in two ways,
$$
\sum_{j \in N} x_j = \sum_{j \in N} \, \sum_{i \in S_j(\bfx)} x_i = \sum_{j \in N} \, \sum_{i \in S_j(\bfx)} j = \sum_{j \in N} j \, |S_j(\bfx)|
= \sum_{j \in N} j \, x'_j .
$$
\end{proof}

\begin{lemma}\label{lemma:g}
Let $\bfx \in N^N$. The following statements hold:\\
\begin{tabular}{rl}
$(i)$ & $s(\bfx) \in N^N$ if and only if $x_0,\ldots,x_{n-1}$ are not all equal.\\
$(ii)$ & If $s(\bfx) \in N^N$ then $s(s(\bfx)) \in N^N$ if and only if $x_0,\ldots,x_{n-1}$ are not all distinct.\\
$(iii)$ & If $s(\bfx),s(s(\bfx)) \in N^N$ then $s(s(s(\bfx))) \in N^N$ if and only if $n \geqslant 4$.
\end{tabular}
\end{lemma}

\begin{proof}
$(i)$ Easy.

$(ii)$ Setting $\bfx' := s(\bfx)$ and $\bfx'' := s(\bfx')$, we have
\begin{eqnarray*}
\bfx'' \in N^N
& \Leftrightarrow & x'_0,\ldots,x'_{n-1} \; \mbox{are not all equal} \qquad \mbox{(by $(i)$)}\\
& \Leftrightarrow & \bfx' \neq (1,\ldots,1) \qquad \mbox{(by Eq.\ (\ref{eq:sum/xj}))} \\
& \Leftrightarrow & \{x_0,\ldots,x_{n-1}\} \neq N.
\end{eqnarray*}

$(iii)$ By $(i)$ and $(ii)$, the numbers $x_0,\ldots,x_{n-1}$ are neither all equal nor all distinct, and hence $n \geqslant 3$. Now set $\bfx'
:= s(\bfx)$, $\bfx'' := s(\bfx')$, and $\bfx''' := s(\bfx'')$. By $(ii)$, we have
$$
\bfx''' \in N^N \enskip \Leftrightarrow \enskip \{x'_0,\ldots,x'_{n-1}\} \neq N.
$$
However we have
\begin{eqnarray*}
\{x'_0,\ldots,x'_{n-1}\} = N
& \Rightarrow & \sum_{j \in N} x'_j = \sum_{j \in N} j \\
& \Rightarrow & n = \frac{n(n-1)}2 \qquad \mbox{(by Eq.\ (\ref{eq:sum/xj}))} \\
& \Rightarrow & n = 3
\end{eqnarray*}
and
$$
n = 3 \enskip \Rightarrow \enskip \{x'_0,x'_1,x'_2\} = \{0,1,2\} = N.
$$
Thus Lemma~\ref{lemma:g} is proved.
\end{proof}

Let $\N$ denote the set of all $n$-chains whose components are neither all equal nor all distinct. One can readily see that $|\N| = n^n - n! -
n$. Moreover, we have the following result, which immediately follows from Lemma~\ref{lemma:g}.

\begin{proposition}\label{prop:sequence}
Let $\bfx \in N^N$. Then all the chains $s(\bfx), s(s(\bfx)), s(s(s(\bfx))),\ldots$ belong to $N^N$ if and only if $\bfx \in \N$ and $n
\geqslant 4$. In that case, all these chains belong to $\N$.
\end{proposition}

From now on we will assume that $n \geqslant 4$. Let $\nit$ denote the set of non-negative integers. According to
Proposition~\ref{prop:sequence} we can construct from any $\bfx \in \N$ an infinite sequence of chains $(\bfx^{(k)})_{k \in \nit}$ in the
following way:
\begin{equation}\label{eq:seq}
\cases{\bfx^{(0)} = \bfx, & \cr \bfx^{(k+1)} = s(\bfx^{(k)}), & $\; k \in \nit$.\cr}
\end{equation}
Since $\N$ is a finite set, this sequence is eventually periodic. That is, there exist $k_0,l \in \nit$ ($l \geqslant 1$) such that
\begin{equation}\label{eq:per}
\bfx^{(k+l)} = \bfx^{(k)} \qquad \forall k \geqslant k_0.
\end{equation}
If the chains $\bfx^{(k)},\ldots,\bfx^{(k+l-1)}$ are distinct and such that $\bfx^{(k+l)} = \bfx^{(k)}$, we say that they form a {\em circuit}
of length $l$. Of course, determining perfect (resp.\ amicable, sociable) chains amounts to identifying all the circuits of length 1 (resp.\ 2,
$\geqslant 3$).

\section{Exhaustive computation of perfect, amicable, and sociable chains}

In the present section we calculate all the perfect, amicable, and sociable chains. These are given in Theorem~\ref{thm:list} below.

Assume that $\bfx^{(k_0)} \in \N$ belongs to a circuit. By Proposition~\ref{prop:sequence}, we have $\bfx^{(k)} \in \N$ for all $k \geqslant
k_0$. Furthermore, by Eq.\ (\ref{eq:sum/xj}) and (\ref{eq:sum/jxj}), we have
\begin{eqnarray}
\sum_{j \in N} x^{(k)}_j & = & n \qquad \forall k \geqslant k_0,\label{eq:sum/xkj} \\
\sum_{j \in N} j \, x^{(k)}_j & = & n \qquad \forall k \geqslant k_0.\label{eq:sum/jxkj}
\end{eqnarray}
These identities imply trivially
\begin{equation}\label{eq:x0}
x^{(k)}_0 = \sum_{j=1}^{n-1} (j-1) \, x^{(k)}_j \qquad \forall k \geqslant k_0.
\end{equation}

Moreover, we have
\begin{equation}\label{eq:x0/1}
x_0^{(k)} \geqslant 1 \qquad \forall k \geqslant k_0.
\end{equation}
Indeed, if $x_0^{(k)} = 0$ for some $k \geqslant k_0$ then, by Eq.~(\ref{eq:x0}), we have $x_2^{(k)}=\cdots =x_{n-1}^{(k)}=0$. By
Eq.~(\ref{eq:sum/xkj}) we then have $x_1^{(k)}=n$, a contradiction.

\begin{theorem}\label{thm:list}
Let $\|$ denote a list, possibly empty, of zeroes.
\par \noindent
The perfect chains are:
\begin{eqnarray}
&& (1,2,1,0) \label{eq:c1}\\
&& (2,0,2,0) \label{eq:c2}\\
&& (2,1,2,0,0) \label{eq:c3}\\
&& (n-4,2,1 \| 1,0,0,0), \qquad n \geqslant 7. \label{eq:c4}
\end{eqnarray}
The pairs of amicable chains are:
\begin{eqnarray}
&& (2,3,0,1,0,0), \quad (3,1,1,1,0,0) \label{eq:c5}\\
&& (n-4,3,0,0 \| 0,1,0,0), \quad (n-3,1,0,1 \| 1,0,0,0) , \qquad n \geqslant 8. \label{eq:c6}
\end{eqnarray}
The unique group of sociable chains is:
\begin{equation}\label{eq:c7}
(3,3,0,0,1,0,0), \quad (4,1,0,2,0,0,0), \quad (4,1,1,0,1,0,0).
\end{equation}
There is no group of more than 3 sociable chains.
\end{theorem}

\begin{proof}
Let $\bfx^{(k_0)} \in \N$ belong to a circuit. Choose $k \geqslant k_0$ such that $x_0^{(k+1)} \leqslant x_0^{(k)}$. Such a $k$ exists for
otherwise $\bfx^{(k_0)}$ would not belong to a circuit.

Set $p := x_0^{(k)}$. By Eq.\ (\ref{eq:x0/1}), we have $1 \leqslant p \leqslant n-1$. Moreover, since $0 \in S_{p}(\bfx^{(k)})$, we have
$$
x_{p}^{(k+1)} = |S_{p}(\bfx^{(k)})| \geqslant 1.
$$

Using Eq.\ (\ref{eq:x0}), we have
$$
x_0^{(k+1)} = \sum_{j=1}^{n-1} (j-1) \, x_j^{(k+1)} \geqslant (p-1) + \sum_{\textstyle{j=1 \atop j \neq p}}^{n-1} (j-1) \, x_j^{(k+1)},
$$
and hence,
\begin{equation}\label{eq:enca}
1 \geqslant 1+x_0^{(k+1)}-p \geqslant \sum_{\textstyle{j=1 \atop j \neq p}}^{n-1} (j-1) \, x_j^{(k+1)} \geqslant 0,
\end{equation}
implying $x_0^{(k+1)} = p$ or $x_0^{(k+1)} = p-1$. We now investigate these two cases separately.
\begin{enumerate}
\item {\it Case $x_0^{(k+1)} = p$}.

By Eq.\ (\ref{eq:enca}), we have
$$
x_j^{(k+1)} = 0 \qquad \forall j \in N \setminus \{0,1,2,p\}.
$$
\begin{enumerate}
\item {\it Case $p = 1$}.

Using Eq.\ (\ref{eq:x0}) and (\ref{eq:sum/xkj}), we obtain $x_2^{(k+1)} = 1$ and $x_1^{(k+1)} = n-2$, so that
$$
\bfx^{(k+1)} = (1,n-2,1 \| )
$$
and $\{x_1^{(k)},\ldots,x_{n-1}^{(k)}\} = \{2,1,\ldots,1,0\}$. By Eq.\ (\ref{eq:sum/jxkj}), we have
$$
n = \sum_{j \in N} j \, x_j^{(k)} \geqslant 2 + \sum_{j=2}^{n-2} j = \frac 12 (n-2)(n-1) + 1,
$$
that is $n = 4$. This leads to the circuit (\ref{eq:c1}).

\item {\it Case $p = 2$}.

Using Eq.\ (\ref{eq:x0}) and (\ref{eq:sum/xkj}), we obtain $x_2^{(k+1)} = 2$ and $x_1^{(k+1)} = n-4$, so that
$$
\bfx^{(k+1)} = (2,n-4,2 \| )
$$
and $\{x_1^{(k)},\ldots,x_{n-1}^{(k)}\} = \{2,1,\ldots,1,0,0\}$. By Eq.\ (\ref{eq:sum/jxkj}), we have
$$
n = \sum_{j \in N} j \, x_j^{(k)} \geqslant 2 + \sum_{j=2}^{n-3} j = \frac 12 (n-3)(n-2) + 1,
$$
that is $n \in \{4,5\}$. This leads to the circuits (\ref{eq:c2}) and (\ref{eq:c3}).

\item {\it Case $p \geqslant 3$}.

By Eq.\ (\ref{eq:x0}), we have
$$
p = x_2^{(k+1)} + (p-1) \, x_{p}^{(k+1)},
$$
which implies $x_2^{(k+1)} = x_{p}^{(k+1)} = 1$. By Eq.\ (\ref{eq:sum/xkj}), we then have $x_1^{(k+1)} = n-p-2$, and hence
$$
\bfx^{(k+1)} = (\underbrace{p,n-p-2,1 \| 1}_{\textstyle{p+1}} \| ),
$$
with $n \geqslant p + 2$.
\begin{enumerate}

\item {\it Case $n = p + 2 \; (\geqslant 5)$}.

We have
$$
\bfx^{(k+1)} = (n-2,0,1 \| 1,0), \quad \bfx^{(k+2)} = (n-3,2,0 \| 1,0).
$$
For $n = 5$, we get the circuit (\ref{eq:c3}). For $n \geqslant 6$, we have
$$
\bfx^{(k+3)} = (n-3,1,1 \| 1,0,0), \quad \bfx^{(k+4)} = (n-4,3,0 \| 1,0,0).
$$
For $n = 6$, $n = 7$ and $n \geqslant 8$, we get the circuits (\ref{eq:c5}), (\ref{eq:c7}) and (\ref{eq:c6}) respectively.

\item {\it Case $n = p + 3 \; (\geqslant 6)$}.

We have
$$
\bfx^{(k+1)} = (n-3,1,1 \| 1,0,0),
$$
which leads to a previous case.

\item {\it Case $n = p + 4 \; (\geqslant 7)$}.

We have
$$
\bfx^{(k+1)} = (n-4,2,1 \| 1,0,0,0),
$$
which leads to the circuit (\ref{eq:c4}).

\item {\it Case $n = p + 5 \; (\geqslant 8)$}.

We have
$$
\bfx^{(k+1)} = (n-5,3,1 \| 1,0,0,0,0).
$$
For $n = 8$, we get the circuit (\ref{eq:c6}). For $n \geqslant 9$, we have
$$
\bfx^{(k+2)} = (n-4,2,0,1 \| 1,0,0,0,0),
$$
retrieving the circuit (\ref{eq:c4}).

\item {\it Case $n = p + r \; (\geqslant 3+r)$, with $r \geqslant 6$}.

We have
$$
\bfx^{(k+1)} = (\underbrace{n-r,r-2,1 \| 1}_{\textstyle{n-r+1}} \|).
$$

If $n-r < r-2$ then
$$
\bfx^{(k+2)} = (\underbrace{\overbrace{n-4,2,0 \| 1}^{\textstyle{n-r+1}} \| 1}_{\textstyle{r-1}} \|),
$$
which leads to a previous case.

If $n-r = r-2$ then
$$
\bfx^{(k+2)} = (\underbrace{n-4,2,0 \| 2}_{\textstyle{n-r+1}} \|), \quad \bfx^{(k+3)} = (n-3,0,2 \| 1,0,0,0),
$$
which leads to a previous case.

If $n-r > r-2$ then
$$
\bfx^{(k+2)} = (\underbrace{\overbrace{n-4,2,0 \| 1}^{\textstyle{r-1}} \| 1}_{\textstyle{n-r+1}} \|),
$$
which leads to a previous case.
\end{enumerate}
\end{enumerate}

\item {\it Case $x_0^{(k+1)} = p-1$}.

By Eq.\ (\ref{eq:enca}), we have
$$
x_j^{(k+1)} = 0 \qquad \forall j \in N \setminus \{0,1,p\},
$$
with $p = x_0^{(k+1)} + 1 \geqslant 2$. Using Eq.\ (\ref{eq:sum/xkj}) and (\ref{eq:x0}), we obtain $x_{p}^{(k+1)} = 1$ and $x_1^{(k+1)} = n-p$,
so that
$$
\bfx^{(k+1)} = (\underbrace{p-1,n-p \| 1}_{\textstyle{p+1}} \| ).
$$
\begin{enumerate}

\item {\it Case $p = 2$}.

We have
$$
\bfx^{(k+1)} = (1,n-2,1 \| ),
$$
that is a case previously encountered.

\item {\it Case $p \geqslant 3$}.
\begin{enumerate}

\item {\it Case $n = p+1 \; (\geqslant 4)$}.

We have
$$
\bfx^{(k+1)} = (n-2,1 \| 1),
$$
which leads to a previous case.

\item {\it Case $n = p+r \; (\geqslant 3+r)$, with $r \geqslant 2$}.

We have
$$
\bfx^{(k+1)} = (\underbrace{n-r-1,r \,\| 1}_{\textstyle{n-r+1}} \| ).
$$

If $n-r-1 < r$ then
$$
\bfx^{(k+2)} = (\underbrace{\overbrace{n-3,1 \| 1}^{\textstyle{n-r}} \| 1}_{\textstyle{r+1}} \| ),
$$
which leads to a previous case.

If $n-r-1 = r$ then
$$
\bfx^{(k+2)} = (\underbrace{n-3,1 \| 2}_{\textstyle{n-r}} \| ),
$$
which leads to a previous case.

If $n-r-1 > r$ then
$$
\bfx^{(k+2)} = (\underbrace{\overbrace{n-3,1 \| 1}^{\textstyle{r+1}} \| 1}_{\textstyle{n-r+1}} \| ),
$$
which leads to a previous case.
\end{enumerate}
\end{enumerate}
\end{enumerate}
Theorem~\ref{thm:list} is now proved.
\end{proof}

\begin{corollary}
Any circuit of length $\geqslant 2$ contains the chain $(n-4,3 \| 1,0,0)$.
\end{corollary}

Before closing this section, we present the following open problem. For any $\bfx \in \N$, we denote by ${\cal C}(\bfx)$ the circuit obtained
from the infinite sequence $(\bfx^{(k)})_{k \in \nit}$. The question then arises of determining the length of the non-periodic part of this
sequence; that is, the number of elements that do not belong to ${\cal C}(\bfx)$:
$$
\Psi(\bfx) := \min \{k \in \nit \mid \bfx^{(k)} \in {\cal C}(\bfx)\}.
$$
Interestingly enough, the following sequence:
$$
\psi(n) := \max_{\bfx \in \N} \, \Psi(\bfx), \qquad n \geqslant 4,
$$
has a rather strange behavior. Its first values (for $4 \leqslant n \leqslant 44$) are: 3, 4, 7, 4, 7, 7, 7, 6, 7, 6, 7, 7, 7, 6, 7, 7, 7, 7, 7,
7, 7, 7, 7, 7, 7, 7, 7, 7, 7, 7, 7, 8, 8, 8, 8, 8, 8, 8, 8, 8, 8.

We conjecture that the elements of this sequence can be arbitrary large; that is, for any $M \geqslant 3$ there exists $n\geqslant 4$ such that
$\psi(n) \geqslant M$.

\section{Range of the counting operator and its iterates}

For any $k \in \nit$, let $s^{(k)}$ denote the $k$th iterate of the operator $s$. It is clear that we have
$$
s^{(k+1)}(\N) \subseteq s^{(k)}(\N), \qquad k \in \nit.
$$

In this final section we intend to describe the subset $s^{(k)}(\N)$ for each $k \in \nit$. The case $k=1$ is dealt with in the next
proposition.

\begin{proposition}\label{prop:s/sum}
We have
$$
s(\N) = \Big\{ \bfx \in \N \, \Big| \, \textstyle{\sum\limits_{j \in N} x_j = n}\Big\}.
$$
\end{proposition}

\begin{proof}
$(\subseteq)$ Follows from Eq.\ (\ref{eq:sum/xj}).\\
$(\supseteq)$ Let $\bfx \in \N$ such that $\sum_{j \in N} x_j = n$. Setting
$$
\bfz :=
(\underbrace{0,\ldots,0}_{\textstyle{x_0}},\underbrace{1,\ldots,1}_{\textstyle{x_1}},\ldots,\underbrace{n-1,\ldots,n-1}_{\textstyle{x_{n-1}}}),
$$
we have $\bfz \in \N$ and $s(\bfz) = \bfx$, and hence $\bfx \in s(\N)$.
\end{proof}

Let the operator $r : s(\N) \to \N$ be defined by
$$
r(\bfx) =
(\underbrace{0,\ldots,0}_{\textstyle{x_0}},\underbrace{1,\ldots,1}_{\textstyle{x_1}},\ldots,\underbrace{n-1,\ldots,n-1}_{\textstyle{x_{n-1}}}).
$$
Let $\Pi_N$ be the set of all the permutations on $N$ and define the operator $q : \N \to \N$ by
$$
q(\bfx) = (x_{\nu(0)},\ldots,x_{\nu(n-1)}),
$$
where $\nu \in \Pi_N$ is such that $x_{\nu(0)} \leqslant \cdots \leqslant x_{\nu(n-1)}$. One can easily see that $s \circ r = \mathrm{id}$ and
$r \circ s = q$, thus showing that $s$ is not invertible.

For any $\pi \in \Pi_N$, we define $r_\pi : s(\N) \to \N$ by
$$
r_\pi(\bfx) = (r(\bfx)_{\pi(0)},\ldots,r(\bfx)_{\pi(n-1)}).
$$
For any $\bfx \in s(\N)$, we clearly have $s^{(-1)}(\bfx) = \{r_\pi(\bfx) \mid \pi \in \Pi_N\}$. Moreover, we have the following result.

\begin{proposition}\label{prop:sk/sum}
For any $k \in \nit$, we have
$$
s^{(k+1)}(\N) = \Big\{\bfx \in s^{(k)}(\N) \, \Big| \, \exists \, \pi_1,\ldots,\pi_k \in \Pi_N : \textstyle{\sum\limits_{j \in N}
(r_{\pi_1}\circ\cdots\circ r_{\pi_k})(\bfx)_j = n} \Big\}.
$$
\end{proposition}

\begin{proof}
We proceed by induction over $k \in \nit$. By Proposition~\ref{prop:s/sum}, the result holds for $k = 0$. Assume that it also holds for $k = 0,\ldots,K-1$, with a given $K \geqslant 1$. We now show that it still holds for $k=K$.\\
$(\subseteq)$ Let $\bfx \in s^{(K+1)}(\N)$. Take $\pi_K \in \Pi_N$ and set $\bfz := r_{\pi_K}(\bfx)$. We have $\bfx = s(\bfz)$ and hence $\bfz
\in s^{(K)}(\N)$. By induction hypothesis, there exist $\pi_1,\ldots,\pi_{K-1} \in \Pi_N$ such that
$$
\sum_{j \in N} (r_{\pi_1}\circ\cdots\circ r_{\pi_K})(\bfx)_j = \sum_{j \in N} (r_{\pi_1}\circ\cdots\circ r_{\pi_{K-1}})(\bfz)_j = n.
$$
\noindent $(\supseteq)$ Let $\bfx \in s^{(K)}(\N)$ and assume that there exist $\pi_1,\ldots,\pi_K \in \Pi_N$ such that
$$
\sum_{j \in N} (r_{\pi_1}\circ\cdots\circ r_{\pi_K})(\bfx)_j = n.
$$
We only have to prove that $\bfx \in \{s(\bfz) \mid \bfz \in s^{(K)}(\N)\}$. Set $\bfz := r_{\pi_K}(\bfx)$. We have $\bfx = s(\bfz)$ and hence
$\bfz \in s^{(K-1)}(\N)$. Moreover, we have
$$
\sum_{j \in N} (r_{\pi_1}\circ\cdots\circ r_{\pi_{K-1}})(\bfz)_j = \sum_{j \in N} (r_{\pi_1}\circ\cdots\circ r_{\pi_K})(\bfx)_j = n,
$$
and hence $\bfz \in s^{(K)}(\N)$ by induction hypothesis.
\end{proof}

The case $k = 2$ is particularly interesting. One can easily see that, for any $\bfx \in \N$ and any $j \in N$, $s^{(2)}(\bfx)_j$ represents the
number of distinct values occuring $j$ times in $\{x_0,\ldots,x_{n-1}\}$. Moreover, we have the following proposition.

\begin{proposition}\label{prop:s2}
We have
$$
s^{(2)}(\N) = \Big\{\bfx \in s(\N) \, \Big| \, \textstyle{\sum\limits_{j\in N} j\, x_j = n}\Big\}.
$$
\end{proposition}

\begin{proof}
For any $\pi \in \Pi_N$, we have
$$
\sum_{j \in N} r_\pi(\bfx)_j = \sum_{j \in N} j \, x_j.
$$
We then conclude by Proposition~\ref{prop:sk/sum}.
\end{proof}

Now, from the identity
$$
\Big|\Big\{\bfx \in \nit^n \, \Big| \, \textstyle{\sum\limits_{j=1}^n j \, x_j = n}\Big\}\Big| = P(n),
$$
where $P(n)$ is the number of unrestricted partitions of the integer $n$ (see e.g.\ \cite{AbSt72}), we can easily show that $|s^{(2)}(\N)| =
P(n)-2$. Similarly, from the well-known identity
$$
\Big|\Big\{\bfx \in \nit^n \, \Big| \, \textstyle{\sum\limits_{j=1}^n x_j = n}\Big\}\Big| = {2n-1 \choose n},
$$
we can readily see that $|s(\N)| = {2n-1 \choose n}-n-1$.

Finally, from the identities $r \circ s = q$ and $r \circ s^{(2)} = q \circ s$, we clearly have
\begin{eqnarray*}
r(s(\N)) &=& \{\bfx \in \N \mid x_0 \leqslant \cdots \leqslant x_{n-1}\},\\
r(s^{(2)}(\N)) &=& \Big\{\bfx \in \N \, \Big| \, \textstyle{\sum\limits_{j \in N} x_j = n} \enskip \mbox{and} \enskip x_0 \leqslant \cdots
\leqslant x_{n-1} \Big\},
\end{eqnarray*}
and, since $r$ is an injection, we have
$$
|r(s(\N))| = |s(\N)| \quad \mbox{and} \quad |r(s^{(2)}(\N))| = |s^{(2)}(\N)|.
$$

\end{document}